\documentclass[12pt,reqno]{article}
\usepackage{mathtools}
\usepackage[usestackEOL]{stackengine}[2013-10-15]
\usepackage[usenames]{color}
\usepackage{amsmath,amsthm,amssymb,tikz,enumerate}
\usetikzlibrary{decorations.pathreplacing}
\usepackage{multicol}

\usepackage{color}
\usepackage{fullpage}
\usepackage{float}
\usepackage{epsf}
\newcommand\tab[1][1cm]{\hspace*{#1}}
\usepackage[colorlinks=true,linkcolor=webgreen,filecolor=webbrown,citecolor=webgreen]{hyperref}
\usepackage{graphicx}

\definecolor{webgreen}{rgb}{0,.5,0}
\definecolor{webbrown}{rgb}{.6,0,0}

\begin{document}

\begin{center}
\end{center}
\theoremstyle{plain}
\newtheorem{theorem}{Theorem}
\newtheorem{conclusion}{Conclusion}
\newtheorem{corollary}{Corollary}
\newtheorem{lemma}{Lemma}
\theoremstyle{definition}
\newtheorem{definition}{Definition}
\begin{center}
\vskip 1cm{\LARGE\bf Investigating First Returns: The Effect of Multicolored Vectors}
\vskip 1cm
\large
Shakuan Frankson and Myka Terry\\ 
Mathematics Department\\
SPIRAL Program at Morgan State University, Baltimore, MD\\
\href{mailto:shakuan.frankson@bison.howard.edu}{\tt shakuan.frankson@bison.howard.edu}\\
\href{myter1@morgan.edu}{\tt myter1@morgan.edu}\\

\end{center}
\vskip .2 in

\begin{abstract}
\noindent By definition, a \textit{first} return is the immediate moment that a path, using vectors in the Cartesian plane, touches the $x$-axis after leaving it previously from a given point; the initial point is often the origin. In this case, using certain diagonal and horizontal vectors while restricting the movements to the first quadrant will cause almost every first return to end at the point $(2n,0)$, where $2n$ counts the equal number of up and down steps in a path. The exception will be explained further in the sections below.
\\
\\
Using the first returns of Catalan, Schr\"{o}der, and Motzkin numbers, which resulted from the lattice paths formed using a combination of diagonal and/or horizontal vectors, we then investigated the effect that coloring these vectors will have on each of their respective generating functions. 
\end{abstract}
\emph{\section{Introduction}}
The main focus of the paper is to analyze the effect of multicoloring certain vectors used to produce the Catalan, Schr\"{o}der, and Motzkin paths. In the process of multicoloring the vectors, we set an arbitrary variable $n$ to represent the number of colors assigned to each vector, where each color will stand as an option for a path to take. Further explanation will be provided in Section 3. In addition, the terms \emph{vector} and \emph{step} will be used interchangeably.\\

Before going into direct detail about the first returns of the Catalan, Schr\"{o}der, and Motzkin paths, a synopsis will also be provided below to give additional background.
\vspace{.3cm}
\subsection{Catalan Numbers}
The original  Dyck  path uses the vectors $(1,1)$, referred to as the up step, and $(1,-1)$, referred to as the down step, remaining above the $x$-axis and extending from $(0,0)$ to $(2n,0)$, where $n$ is the semi-length. The Catalan numbers, denoted $C(x)$, are the values enumerated along the bottom row of the x-axis using those vectors. The following vectors are visually represented in the lattice path as:
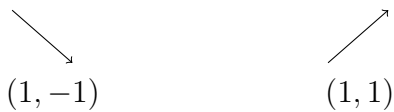
\begin{center}
\begin{tikzpicture}
	\draw[->](3,0) -- (3.8,-0.7);
	\draw[->](7.2,-0.7) -- (8,0);
	
\end{tikzpicture}		

\tab[0.00015cm]$(1,-1)$\tab[3cm]$(1,1)$\\ 
\end{center}

Below is the graph of the first return of the Catalan path using the vectors provided:

\begin{center}
\begin{tikzpicture}[ultra thick]
\draw[->] [color=gray] (0,0)--(0,4);
\draw[->] [color=gray] (0,0)--(10,0);
\draw[->] (0,0)--(1,1) node[midway,above]{\large{$x$} \ };
\draw [nearly transparent, dashed] (1,1) -- (5,1);
\draw[color=teal] (1,1) sin (1.5,2.1) cos (2,1.3) sin (2.5,2.5) cos (3,1.9) sin (3.5,3.5) cos (4,1.85) sin (4.5,2.7) cos (5,1);
\draw[->] (5,1)--(6,0);
\draw[color=teal] (6,0) sin (6.5,1.5) cos (7,0.25) sin (7.5,2.85) cos (8,0.25) sin (8.5,1.85) cos (9,0.4) sin (9.5,1.5) cos (10,0);
\draw[-] (2.2,2.6)--(2.2,2.6) node[midway, above]{\large{$C(x)$}};
\draw[-] (7.4,3)--(7.4,3) node[midway, above]{\large{$C(x)$}};

\end{tikzpicture}

\begin{figure}[h!]
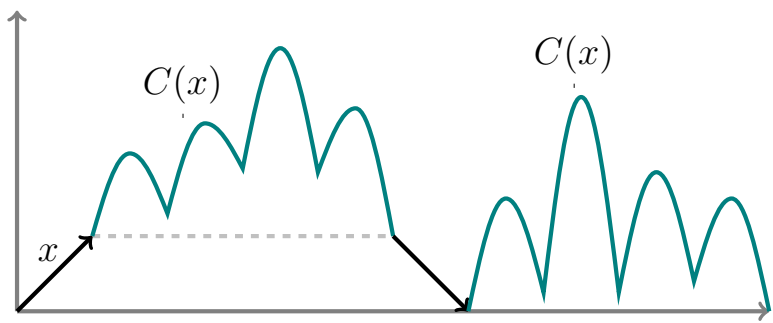

\caption{Example of First Return - $C(x)$}
\end{figure}
\end{center}

To produce the first return for C(x), movement from the origin begins with 2 options: (1) remaining at the origin and not using the provided vectors, or (2) following the $(1,1)$ up step, where we attain an $x$, continuing along a random path of $C(x)$, and returning to the $x$-axis using the $(1,-1)$ down step, ending with an arbitrary infinite path of $C(x)$. Option (1) has the resulting value of $1$ since there is only one way to remain at the origin, and option (2) has a value of $xC(x)\cdot C(x)$. Combining the two possibilities, the first full return is notated as the following: 

$$C(x)=1+xC(x)\cdot C(x) = 1+xC^2(x)$$

The generating function follows from the quadratic formula and yields Equation 1:
\begin{equation}
C(x)=\frac{1-\sqrt{1-4x}}{2x} \label{eq:1}
\end{equation}

\subsection{Schr\"{o}der Numbers}
Schr\"{o}der  numbers follow a similar pattern to the Catalan numbers while also including the level step, $(2,0)$. The path is referred to as a Schr\"{o}der path. There are two distinct types of Schr\"{o}der numbers seen in combinatorics: large Schr\"{o}der numbers, denoted by $S(x)$, and small Schr\"{o}der numbers, denoted by $s(x)$. The following vectors are visually represented in the lattice path as:
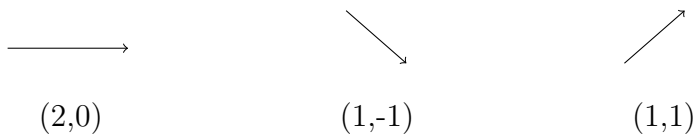
\begin{center}
\begin{tikzpicture}
	\draw[->](0,-0.5) -- (1.6,-0.5);
	\draw[->](4.5,0) -- (5.3,-0.7);
	\draw[->](8.2,-0.7) -- (9,0);
	
\end{tikzpicture}		
\end{center}

\tab[3.55cm](2,0)\tab[3.15cm](1,-1)\tab[2.9cm](1,1)\\

These vectors create unique Schr\"{o}der  paths for the large Schr\"{o}der and small Schr\"{o}der numbers, respectively. The sole distinction between them is $s(x)$ does not consist of the level step, $(2,0)$, along the x-axis, while $S(x)$ does.\\

Below is the graph of the first return of the large Schr\"{o}der path using the vectors provided:

\begin{center}
\begin{tikzpicture}[ultra thick]
\draw[->] [color=gray] (0,0)--(0,4);
\draw[->] [color=gray] (0,0)--(10,0);
\draw[->] (0,0)--(1,1) node[midway,above]{\large{$x$} \ \ \ };
\draw [nearly transparent, dashed] (1,1) -- (5,1);
\draw[->] (0,0)--(2,0);
\draw[color=teal] (1,1) sin (1.5,2.1) cos (2,1.3) sin (2.5,2.5) cos (3,1.9) sin (3.5,3.5) cos (4,1.85) sin (4.5,2.7) cos (5,1);
\draw[->] (5,1)--(6,0);
\draw[color=teal] (6,0) sin (6.5,1.5) cos (7,0.25) sin (7.5,2.85) cos (8,0.25) sin (8.5,1.85) cos (9,0.4) sin (9.5,1.5) cos (10,0);
\draw[double, dotted, color=purple] (2,0) sin (2.5,1.4) cos (3,0.5) sin (3.5,1) cos (4,0) sin (4.5,0.9) cos (5,0.3) sin (5.5,1.75) cos (6,1.15) sin (6.5,2.7) cos (7,.35) sin (7.5,1.35) cos (8,0.5) sin (8.5,0.1) cos (9,0.6) sin (9.5,2.5) cos (10,0);
\draw[-] (2.2,2.6)--(2.2,2.6) node[midway, above]{\large{$S(x)$}};
\draw[-] (7.4,3)--(7.4,3) node[midway, above]{\large{$S(x)$}};
\draw[-] (3.5,1)--(3.5,1) node[midway, above]{\large{$S(x)$}};
\draw[->] (0,0)--(2,0) node[midway,below]{\large{$x$}};

\end{tikzpicture}

\begin{figure}[h!]
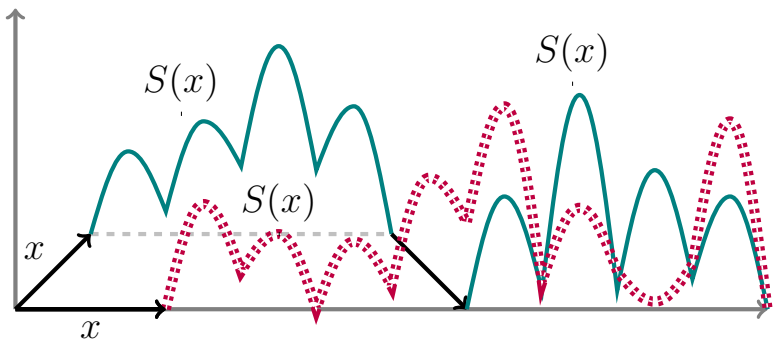

\caption{Example of First Return - $S(x)$}
\end{figure}
\end{center}

To produce the first return for S(x), movement from the origin begins with 3 options: (1) remaining at the origin and not using the provided vectors, (2) following the $(2,0)$ level step, where we attain an $x$, and all of the possibilities of $S(x)$, or (3) following the $(1,1)$ up step, where we attain an $x$, continuing along a random path of $S(x)$, and returning to the $x$-axis using the $(1,-1)$ down step, ending with an arbitrary infinite path of $S(x)$. Option (1) has the resulting value of $1$ since there is only one way to remain at the origin, option (2) has the value of $xS(x)$, and option (3) has a value of $xS(x)\cdot S(x)$. Combining the 3 possibilities, the first full return is notated as the following:
$$S(x)=1+xS(x)+xS(x)\cdot S(x)$$

The generating function follows from the quadratic formula and yields Equation 2:

\begin{equation}
S(x)=\frac{1-x-\sqrt{1-6x+x^2}}{2x} \label{eq:2}
\end{equation}\\

Since the small Schr\"{o}der path does not include the level step, $(2,0)$, the first return appears similar to the first return of the Catalan path. Below is the graph for the first return of small Schr\"{o}der paths:

\begin{center}
\begin{tikzpicture}[ultra thick]
\draw[->] [color=gray] (0,0)--(0,4);
\draw[->] [color=gray] (0,0)--(10,0);
\draw[->] (0,0)--(1,1) node[midway,above]{\large{$x$} \ \ };
\draw [nearly transparent, dashed] (1,1) -- (5,1);
\draw[color=teal] (1,1) sin (1.5,2.1) cos (2,1.3) sin (2.5,2.5) cos (3,1.9) sin (3.5,3.5) cos (4,1.85) sin (4.5,2.7) cos (5,1);
\draw[->] (5,1)--(6,0);
\draw[color=orange] (6,0) sin (6.5,1.5) cos (7,0.25) sin (7.5,2.85) cos (8,0.25) sin (8.5,1.85) cos (9,0.4) sin (9.5,1.5) cos (10,0);
\draw[-] (2.2,2.6)--(2.2,2.6) node[midway, above]{\textbf{\large{$S(x)$}}};
\draw[-] (7.4,3)--(7.4,3) node[midway, above]{{\textbf{\large{$s(x)$}}}};

\end{tikzpicture}

\begin{figure}[h!]
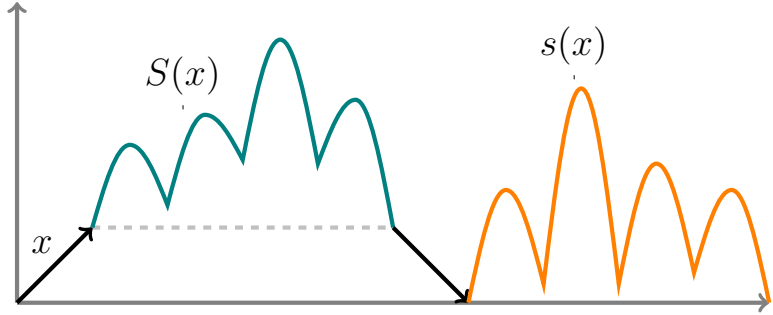

\caption{Example of First Return - $s(x)$}
\end{figure}
\end{center}

The small Schr\"{o}der numbers have a slightly different pattern than the large Schr\"{o}der numbers because they do not have a level step on the $x$-axis, as it affects the first return and therefore its generating function.\\

To produce the first return for s(x), movement from the origin begins with 2 options: (1) remaining at the origin and not using the provided vectors, and (2) following the $(1,1)$ up step, where we attain an $x$, continuing along a random path of $S(x)$, and returning to the $x$-axis using the $(1,-1)$ down step, ending with an arbitrary infinite path of $s(x)$. Option (1) has the resulting value of $1$ since there is only one way to remain at the origin, and option (2) has a value of $xS(x)\cdot s(x)$. Combining these possibilities, the first full return of $s(x)$ is notated as the following:
$$s(x)=1+xS(x)\cdot s(x)$$

Substituting the generating function of $S(x)$ into the equation above and simplifying, we get the generating function for small Schr\"{o}der numbers, denoted as Equation 3:

\begin{equation}
s(x)=\frac{1+x-\sqrt{1-6x+x^2}}{4x} \label{eq:3}
\end{equation}

\subsection{Motzkin Numbers}
Motzkin paths, denoted $M(x)$, are similar in nature to the previous paths, but the  Motzkin  numbers have the smaller level  step  of  $(1,0)$ and remain above the $x$ -axis from $(0,0)$ to $(n,0)$. The following vectors are visually represented in the lattice path as:

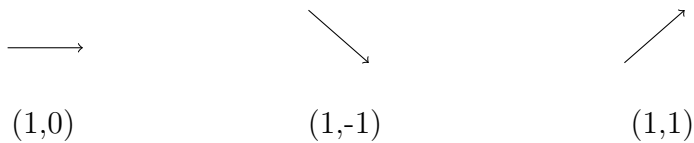
\begin{center}

\begin{tikzpicture}

\draw[->] (0,-.5) -- (1,-.5);

\draw [->](4,0) -- (4.8,-0.7);

\draw [->](8.2,-0.7) -- (9, 0);
\end{tikzpicture}
\end{center}
\tab[3.8cm](1,0)\tab[3.1cm](1,-1)\tab[3.3cm](1,1)\\

Because the length of the horizontal vector is only half that of the $(2,0)$ level step found in the Schr\"{o}der numbers, there are twice as many of these vectors in Motzkin paths.

\begin{center}
\begin{tikzpicture}[ultra thick]
\draw[->] [color=gray] (0,0)--(0,4);
\draw[->] [color=gray] (0,0)--(10,0);
\draw[->] (0,0)--(1,0) node[midway,below]{\large{$x$}};
\draw[->] (5,1)--(6,0) node[midway,above]{\ \ \large{$x$}};
\draw [thick, nearly transparent, dashed] (1,1) -- (5,1);
\draw[->] (0,0)--(1,0);
\draw[color=teal] (1,1) sin (1.5,2.1) cos (2,1.3) sin (2.5,2.5) cos (3,1.9) sin (3.5,3.5) cos (4,1.85) sin (4.5,2.7) cos (5,1);
\draw[->] (5,1)--(6,0);
\draw[color=teal] (6,0) sin (6.5,1.5) cos (7,0.25) sin (7.5,2.85) cos (8,0.25) sin (8.5,1.85) cos (9,0.4) sin (9.5,1.5) cos (10,0);
\draw[double, dotted, color=green] (1,0) sin (1.5,1) cos (2,0) sin (2.5,1.4) cos (3,0.5) sin (3.5,1) cos (4,0) sin (4.5,0.9) cos (5,0.3) sin (5.5,1.75) cos (6,1.15) sin (6.5,2.7) cos (7,.35) sin (7.5,1.35) cos (8,0.5) sin (8.5,0.1) cos (9,0.6) sin (9.5,2.5) cos (10,0);
\draw[-] (2.2,2.6)--(2.2,2.6) node[midway, above]{\large{$M(x)$}};
\draw[-] (7.4,3)--(7.4,3) node[midway, above]{\large{$M(x)$}};
\draw[-] (3.5,1)--(3.5,1) node[midway, above]{\large{$M(x)$}};
\draw [->] (0,0)--(1,1) node[midway,above]{\large{$x$} \ \ };

\end{tikzpicture}

\begin{figure}[h!]
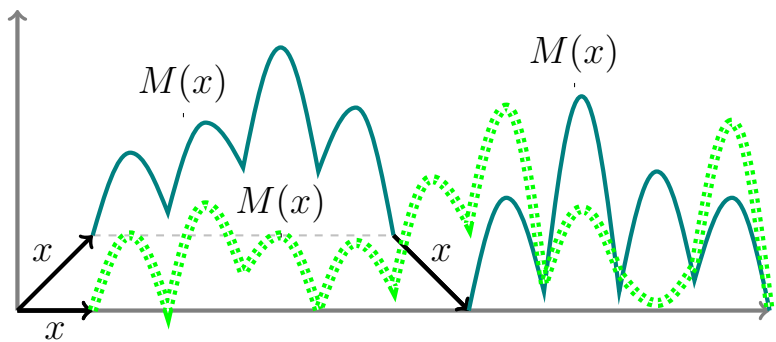

\caption{Example of First Return - $M(x)$}
\end{figure}
\end{center}

To produce the first return for M(x), movement from the origin begins with 3 options: (1) remaining at the origin and not using the provided vectors, (2) following the $(1,0)$ level step, where we attain an $x$, and all of the possibilities of $M(x)$, or (3) following the $(1,1)$ up step, where we attain an $x$, continuing along a random path of $M(x)$, and returning to the $x$-axis using the $(1,-1)$ down step, where we attain another $x$, ending with an arbitrary infinite path of $M(x)$. Option (1) has the resulting value of $1$ since there is only one way to remain at the origin, option (2) has the value of $xM(x)$, and option (3) has a value of $x^2M(x)\cdot M(x)$. 
Combining these three possibilities, the movements of $M(x)$ based on the first movement yield the function: 
$$M(x)=1+xM(x)+x^2M(x)\cdot M(x) = 1+xM(x)+x^2M^2(x)$$

$M(x)$ can be subtracted from both sides, creating a function that can be solved using the quadratic formula, yielding Equation 4:

\begin{equation}
M(x)=\frac{1-x-\sqrt{1-2x-3x^2}}{2x^2} \label{eq:4}
\end{equation}

\section{Effect of Multicolored Vectors on First Returns}
The concept of multicoloring vectors was inspired by an article called \emph{Zebra Trees} written by Davenport, Shapiro, and Woodson, see \cite{DLL}. Zebra trees are ordered trees that have a choice of $2$ colors, black and white, at each edge of odd height, and the coloring of each edge affects the generating function. Expanding off of their research, we decided to investigate the effect of assigning $m$ colors to the down steps and $n$ colors to the level steps for each of the lattice paths mentioned in Section 1. Catalan, Schr\"{o}der, and Motzkin numbers have generating functions that also represent the presence of $1$ color at each vector, which also counts as one choice for a path to take during a return.\\

By finding the generalized generating functions for $m$-colored down steps, $n$-colored level steps, or a combination of both, we will find an efficient method for quickly calculating the values produced by the lattice paths without counting out each individual path one-by-one. The only path that will not be accounted for is the small Schr\"{o}der path since its generalized generating function was more complex to find.

\subsection{Catalan Numbers}
According to Section 1, the first return for a Catalan path is represented by:
\begin{align*}
    C(x) = 1 + xC^2(x)
\end{align*}
where the $x$ is attained from the up step, and $1$ (and hence $1$ $\cdot$ $x$ = $x$) represents that there is $1$ color assigned to the down step.\\

Since the Catalan path only uses the $(1,1)$ up vector and the $(1,-1)$ down vector, then we are only able to find the generating function for $C_{m}(x)$ or $C_{m}$, where $C_{m}$ denotes the Catalan path with $m$-colored $(1,-1)$ down steps. The path for $C_m$ will be represented by:
\begin{align}
    C_{m} = 1+mx\ C^2_{m}
\end{align}

\begin{center}
\begin{tikzpicture}[ultra thick]
\draw[->] [color=gray] (0,0)--(0,4);
\draw[->] [color=gray] (0,0)--(10,0);
\draw[->] (0,0)--(1,1) node[midway,above]{\large{$x$} \ };
\draw [nearly transparent, dashed] (1,1) -- (5,1);
\draw[color=teal] (1,1) sin (1.5,2.1) cos (2,1.3) sin (2.5,2.5) cos (3,1.9) sin (3.5,3.5) cos (4,1.85) sin (4.5,2.7) cos (5,1);
\draw[->, red] (5.1,1)--(6.1,0);
\draw[->, green] (5,1)--(6,0);
\draw[->, yellow] (4.9,1)--(5.9,0);
\draw[->, blue] (4.8,1)--(5.8,0);
\draw[-] (5.65,0.55)--(5.65,0.55) node[midway, above]{\small{$m$}};

\draw[color=teal] (6,0) sin (6.5,1.5) cos (7,0.25) sin (7.5,2.85) cos (8,0.25) sin (8.5,1.85) cos (9,0.4) sin (9.5,1.5) cos (10,0);
\draw[-] (2.2,2.6)--(2.2,2.6) node[midway, above]{\large{$C_m(x)$}};
\draw[-] (7.4,3)--(7.4,3) node[midway, above]{\large{$C_m(x)$}};

\end{tikzpicture}

\begin{figure}[h!]
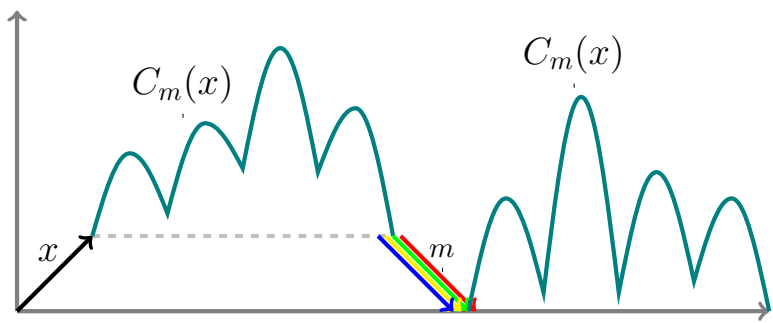

\caption{Example of First Return - $C_m(x)$}
\end{figure}
\end{center}

Now, the $x$ is attained from the up step, and $m$ (and hence $m$ $\cdot$ $x$ = $mx$) represents that there are $m$ colors assigned to the down step.\\

After using the quadratic formula, the generating function yielded for $C_m$ is:

\begin{align}
    C_m = \frac{1-\sqrt{1-4mx}}{2mx}
\end{align}

The generating function for the Catalan path is stated below for comparison purposes:
\begin{align}
    C(x)=\frac{1-\sqrt{1-4x}}{2x} \tag{\ref{eq:1}}
\end{align}
\subsection{Large Schr\"{o}der Numbers}
According to Section 1, the first return for a large Schr\"{o}der path is represented by:
\begin{align*}
    S(x)=1+xS(x)+xS^2(x)
\end{align*}

Large Schr\"{o}der paths use the $(1,1)$ up vector, $(1,-1)$ down vector, and $(2,0)$ horizontal vector. Hence, we are able to find the generating functions for $S_{m}(x)$ or $S_{m}$, where $S_{m}$ denotes the large Schr\"{o}der path with $m$-colored $(1,-1)$ down steps, and $S_{n}(x)$ or $S_{n}$, where $S_{n}$ denotes the large Schr\"{o}der path with $n$-colored $(2,0)$ horizontal steps. The  path for $S_m$ will be represented by:
\begin{align}
    S_{m} = 1+xS_{m}+mx\ S^2_{m} 
\end{align}

\begin{center}
\begin{tikzpicture}[ultra thick]
\draw[->] [color=gray] (0,0)--(0,4);
\draw[->] [color=gray] (0,0)--(10,0);
\draw[->] (0,0)--(1,1) node[midway,above]{\large{$x$} \ \ \ };
\draw [nearly transparent, dashed] (1,1) -- (5,1);
\draw[->] (0,0)--(2,0);
\draw[color=teal] (1,1) sin (1.5,2.1) cos (2,1.3) sin (2.5,2.5) cos (3,1.9) sin (3.5,3.5) cos (4,1.85) sin (4.5,2.7) cos (5,1);
\draw[->, red] (5.1,1)--(6.1,0);
\draw[->, green] (5,1)--(6,0);
\draw[->, yellow] (4.9,1)--(5.9,0);
\draw[->, blue] (4.8,1)--(5.8,0);
\draw[-] (5.65,0.55)--(5.65,0.55) node[midway, above]{\small{$m$}};
\draw[color=teal] (6,0) sin (6.5,1.5) cos (7,0.25) sin (7.5,2.85) cos (8,0.25) sin (8.5,1.85) cos (9,0.4) sin (9.5,1.5) cos (10,0);
\draw[double, dotted, color=purple] (2,0) sin (2.5,1.4) cos (3,0.5) sin (3.5,1) cos (4,0.1) sin (4.5,0.9) cos (5,0.3) sin (5.5,1.75) cos (6,1.15) sin (6.5,2.7) cos (7,.35) sin (7.5,1.35) cos (8,0.5) sin (8.5,0.1) cos (9,0.6) sin (9.5,2.5) cos (10,0);
\draw[-] (2.2,2.6)--(2.2,2.6) node[midway, above]{\large{$S_m(x)$}};
\draw[-] (7.4,3)--(7.4,3) node[midway, above]{\large{$S_m(x)$}};
\draw[-] (3.5,1)--(3.5,1) node[midway, above]{\large{$S_m(x)$}};
\draw[->] (0,0)--(2,0) node[midway,below]{\large{$x$}};

\end{tikzpicture}
\begin{figure}[h!]
\caption{Example of First Return - $S_m(x)$}
\end{figure}
\end{center}

and the path for $S_n$ will be represented by:
\begin{align}
    S_n = 1 + nxS_n+x\ S^2_n
\end{align}

\begin{center}
\begin{tikzpicture}[ultra thick]
\draw[->] [color=gray] (0,0)--(0,4);
\draw[->] [color=gray] (0,0)--(10,0);
\draw[->] (0,0)--(1,1) node[midway,above]{\large{$x$} \ \ \ };
\draw [nearly transparent, dashed] (1,1) -- (5,1);
\draw[color=teal] (1,1) sin (1.5,2.1) cos (2,1.3) sin (2.5,2.5) cos (3,1.9) sin (3.5,3.5) cos (4,1.85) sin (4.5,2.7) cos (5,1);
\draw[->, red] (0,0.09)--(2,0.09);
\draw[->, green] (0,0.025)--(2,0.025);
\draw[->, yellow] (0,-0.05)--(2,-0.05);
\draw[->, blue] (0,-0.12)--(2,-0.12);
\draw[-] (1,-0.15)--(1,-0.15) node[midway, below]{\small{$nx$}};
\draw[->] (5,1)--(6,0);
\draw[color=teal] (6,0) sin (6.5,1.5) cos (7,0.25) sin (7.5,2.85) cos (8,0.25) sin (8.5,1.85) cos (9,0.4) sin (9.5,1.5) cos (10,0);
\draw[double, dotted, color=purple] (2,0) sin (2.5,1.4) cos (3,0.5) sin (3.5,1) cos (4,0.1) sin (4.5,0.9) cos (5,0.3) sin (5.5,1.75) cos (6,1.15) sin (6.5,2.7) cos (7,.35) sin (7.5,1.35) cos (8,0.5) sin (8.5,0.1) cos (9,0.6) sin (9.5,2.5) cos (10,0);
\draw[-] (2.2,2.6)--(2.2,2.6) node[midway, above]{\large{$S_n(x)$}};
\draw[-] (7.4,3)--(7.4,3) node[midway, above]{\large{$S_n(x)$}};
\draw[-] (3.5,1)--(3.5,1) node[midway, above]{\large{$S_n(x)$}};

\end{tikzpicture}

\begin{figure}[h!]
\caption{Example of First Return - $S_n(x)$}
\end{figure}
\end{center}

The path $S_{m,n}$, which denotes the large Schr\"{o}der path with both $m$-colored down steps and $n$-colored level steps, will be represented by:
\begin{align}
    S_{m,n} = 1+nx\ S_{m,n}+mx\ S^2_{m,n}
\end{align}

\begin{center}
\begin{tikzpicture}[ultra thick]
\draw[->] [color=gray] (0,0)--(0,4);
\draw[->] [color=gray] (0,0)--(10,0);
\draw[->] (0,0)--(1,1) node[midway,above]{\large{$x$} \ \ \ };
\draw [nearly transparent, dashed] (1,1) -- (5,1);
\draw[color=teal] (1,1) sin (1.5,2.1) cos (2,1.3) sin (2.5,2.5) cos (3,1.9) sin (3.5,3.5) cos (4,1.85) sin (4.5,2.7) cos (5,1);
\draw[->, red] (0,0.09)--(2,0.09);
\draw[->, green] (0,0.025)--(2,0.025);
\draw[->, yellow] (0,-0.05)--(2,-0.05);
\draw[->, blue] (0,-0.12)--(2,-0.12);
\draw[-] (1,-0.15)--(1,-0.15) node[midway, below]{\small{$nx$}};
\draw[->, red] (5.1,1)--(6.1,0);
\draw[->, green] (5,1)--(6,0);
\draw[->, yellow] (4.9,1)--(5.9,0);
\draw[->, blue] (4.8,1)--(5.8,0);
\draw[-] (5.65,0.55)--(5.65,0.55) node[midway, above]{\small{$m$}};
\draw[color=teal] (6,0) sin (6.5,1.5) cos (7,0.25) sin (7.5,2.85) cos (8,0.25) sin (8.5,1.85) cos (9,0.4) sin (9.5,1.5) cos (10,0);
\draw[double, dotted, color=purple] (2,0) sin (2.5,1.4) cos (3,0.5) sin (3.5,1) cos (4,0.1) sin (4.5,0.9) cos (5,0.3) sin (5.5,1.75) cos (6,1.15) sin (6.5,2.7) cos (7,.35) sin (7.5,1.35) cos (8,0.5) sin (8.5,0.1) cos (9,0.6) sin (9.5,2.5) cos (10,0);
\draw[-] (2.2,2.6)--(2.2,2.6) node[midway, above]{\large{$S_{m,n}(x)$}};
\draw[-] (7.4,3)--(7.4,3) node[midway, above]{\large{$S_{m,n}(x)$}};
\draw[-] (3.5,1)--(3.5,1) node[midway, above]{\large{$S_{m,n}(x)$}};

\end{tikzpicture}

\begin{figure}[h!]
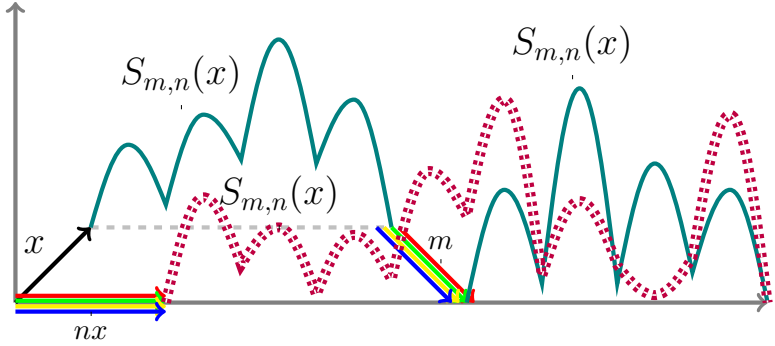

\caption{Example of First Return - $S_{m,n}(x)$}
\end{figure}
\end{center}

After using the quadratic formula, the respective generating functions yielded for $S_m$, $S_n$, and $S_{m,n}$ are:

\begin{align}
    S_m = \frac{1-x-\sqrt{1-2x(1+2m)+x^2}}{2mx}
\end{align}

\begin{align}
    S_n = \frac{1-nx-\sqrt{1-2x(n+2)+x^2n^2}}{2x}
\end{align}

\begin{align}
    S_{m,n} = \frac{1-nx-\sqrt{1-2x(n+2m)+x^2n^2}}{2mx}
\end{align}

The generating function for the large Schr\"{o}der path is stated below for comparison purposes:
\begin{align}
    S(x)=\frac{1-x-\sqrt{1-6x+x^2}}{2x} \tag{\ref{eq:2}}
\end{align}

\subsection{Small Schr\"{o}der Numbers}
According to Section 1, the first return for a small Schr\"{o}der path is represented by:
\begin{align*}
    s(x)=1+xs(x)\cdot S(x)
\end{align*}

The small Schr\"{o}der path uses the $(1,1)$ up vector, $(1,-1)$ down vector, and the $(2,0)$ horizontal vector, but there are no level steps on the $x$-axis. We were only able to find the generating function for $s_{m}(x)$ or $s_{m}$, where $s_{m}$ denotes the small Schr\"{o}der path with $m$-colored $(1,-1)$ down steps. The path for $s_m$ will be represented by:
\begin{align}
    s_{m} = 1+mx\ S_m\cdot s_m
\end{align}

\begin{center}
\begin{tikzpicture}[ultra thick]
\draw[->] [color=gray] (0,0)--(0,4);
 \draw[->] [color=gray] (0,0)--(10,0);
 \draw[->] (0,0)--(1,1) node[midway,above]{\large{$x$}  \ \ };
 \draw [dashed] (1,1) -- (5,1);
 \draw[color=teal] (1,1) sin (1.5,2.1) cos (2,1.3) sin (2.5,2.5) cos (3,1.9) sin (3.5,3.5) cos (4,1.85) sin (4.5,2.7) cos (5,1);
 \draw[->, red] (5.1,1)--(6.1,0);
 \draw[->, green] (5,1)--(6,0);
 \draw[->, yellow] (4.9,1)--(5.9,0);
 \draw[->, blue] (4.8,1)--(5.8,0);
 \draw[-] (5.65,0.55)--(5.65,0.55) node[midway, above]{\small{$m$}};
 \draw[color=orange] (6,0) sin (6.5,1.5) cos (7,0.25) sin (7.5,2.85) cos (8,0.25) sin (8.5,1.85) cos (9,0.4) sin (9.5,1.5) cos (10,0);
 \draw[-] (2.2,2.6)--(2.2,2.6) node[midway, above]{\textbf{\large{$S_m(x)$}}};
 \draw[-] (7.4,3)--(7.4,3) node[midway, above]{{\textbf{\large{$s_m(x)$}}}};

 \end{tikzpicture}

 \begin{figure}[h!]
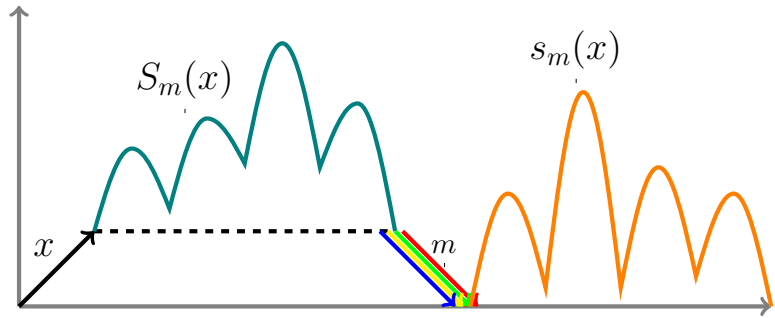

 \caption{Example of First Return - $s_m(x)$}
 \end{figure}
 \end{center}

After substituting in the generating function for $S_m$, the generating function yielded for $s_m$ is:

 \begin{align}
     s_m = \frac{1+x-\sqrt{1-2x(1+2m)+x^2}}{2x(1+m)}
 \end{align}

 The generating function for the small Schr\"{o}der path is stated below for comparison purposes:
 \begin{align}
     s(x)=\frac{1+x-\sqrt{1-6x+x^2}}{4x} \tag{\ref{eq:3}}
 \end{align}

\subsection{Motzkin Numbers}
According to Section 1, the first return for a Motzkin path is represented by:
\begin{align*}
   M(x)=1+xM(x)+x^2M^2(x)
\end{align*}

Motzkin paths use the $(1,1)$ up vector, $(1,-1)$ down vector, and $(1,0)$ horizontal vector. Hence, we are able to find the generating functions for $M_{m}(x)$ or $M_{m}$, where $M_{m}$ denotes the large Schr\"{o}der path with $m$-colored $(1,-1)$ down steps, and $M_{n}(x)$ or $M_{n}$, where $M_{n}$ denotes the large Schr\"{o}der path with $n$-colored $(1,0)$ horizontal steps. The  path for $M_m$ will be represented by:
\begin{align}
    M_{m} = 1+xM_{m}+mx^2\ M^2_{m} 
\end{align}

\begin{center}
\begin{tikzpicture}[ultra thick]
\draw[->] [color=gray] (0,0)--(0,4);
\draw[->] [color=gray] (0,0)--(10,0);
\draw[->] (0,0)--(1,0) node[midway,below]{\large{$x$}};
\draw[->] (5,1)--(6,0);
\draw [thick, nearly transparent, dashed] (1,1) -- (5,1);
\draw[->] (0,0)--(1,0);
\draw[color=teal] (1,1) sin (1.5,2.1) cos (2,1.3) sin (2.5,2.5) cos (3,1.9) sin (3.5,3.5) cos (4,1.85) sin (4.5,2.7) cos (5,1);
\draw[->, red] (5.1,1)--(6.1,0);
\draw[->, green] (5,1)--(6,0);
\draw[->, yellow] (4.9,1)--(5.9,0);
\draw[->, blue] (4.8,1)--(5.8,0);
\draw[-] (5.65,0.55)--(5.65,0.55) node[midway, above]{\small{$mx$}};
\draw[color=teal] (6,0) sin (6.5,1.5) cos (7,0.25) sin (7.5,2.85) cos (8,0.25) sin (8.5,1.85) cos (9,0.4) sin (9.5,1.5) cos (10,0);
\draw[double, dotted, color=green] (1,0) sin (1.5,1) cos (2,0.2) sin (2.5,1.4) cos (3,0.5) sin (3.5,1) cos (4,0.1) sin (4.5,0.9) cos (5,0.3) sin (5.5,1.75) cos (6,1.15) sin (6.5,2.7) cos (7,.35) sin (7.5,1.35) cos (8,0.5) sin (8.5,0.1) cos (9,0.6) sin (9.5,2.5) cos (10,0);
\draw[-] (2.2,2.6)--(2.2,2.6) node[midway, above]{\large{$M_m(x)$}};
\draw[-] (7.4,3)--(7.4,3) node[midway, above]{\large{$M_m(x)$}};
\draw[-] (3.5,1)--(3.5,1) node[midway, above]{\large{$M_m(x)$}};
\draw [->] (0,0)--(1,1) node[midway,above]{\large{$x$} \ \ };

\end{tikzpicture}
\begin{figure}[h!]
\caption{Example of First Return - $M_m(x)$}
\end{figure}
\end{center}

and the path for $M_n$ will be represented by:
\begin{align}
    M_n = 1 + nxM_n+x^2\ M^2_n
\end{align}

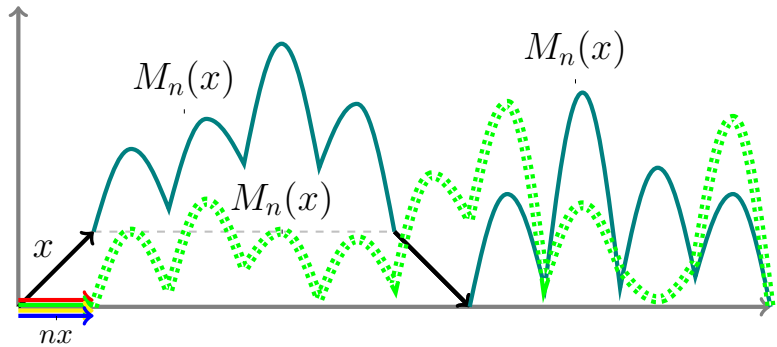
\begin{center}
\begin{tikzpicture}[ultra thick]
\draw[->] [color=gray] (0,0)--(0,4);
\draw[->] [color=gray] (0,0)--(10,0);
\draw[->] (5,1)--(6,0);
\draw [thick, nearly transparent, dashed] (1,1) -- (5,1);
\draw[color=teal] (1,1) sin (1.5,2.1) cos (2,1.3) sin (2.5,2.5) cos (3,1.9) sin (3.5,3.5) cos (4,1.85) sin (4.5,2.7) cos (5,1);
\draw [->] (0,0)--(1,1) node[midway,above]{\large{$x$} \ \ };
\draw[->] (5,1)--(6,0);
\draw[->, red] (0,0.09)--(1,0.09);
\draw[->, green] (0,0.025)--(1,0.025);
\draw[->, yellow] (0,-0.05)--(1,-0.05);
\draw[->, blue] (0,-0.12)--(1,-0.12);
\draw[-] (0.5,-0.15)--(0.5,-0.15) node[midway, below]{\small{$nx$}};
\draw[color=teal] (6,0) sin (6.5,1.5) cos (7,0.25) sin (7.5,2.85) cos (8,0.25) sin (8.5,1.85) cos (9,0.4) sin (9.5,1.5) cos (10,0);
\draw[double, dotted, color=green] (1,0) sin (1.5,1) cos (2,0.2) sin (2.5,1.4) cos (3,0.5) sin (3.5,1) cos (4,0.1) sin (4.5,0.9) cos (5,0.3) sin (5.5,1.75) cos (6,1.15) sin (6.5,2.7) cos (7,.35) sin (7.5,1.35) cos (8,0.5) sin (8.5,0.1) cos (9,0.6) sin (9.5,2.5) cos (10,0);
\draw[-] (2.2,2.6)--(2.2,2.6) node[midway, above]{\large{$M_n(x)$}};
\draw[-] (7.4,3)--(7.4,3) node[midway, above]{\large{$M_n(x)$}};
\draw[-] (3.5,1)--(3.5,1) node[midway, above]{\large{$M_n(x)$}};
\end{tikzpicture}
\begin{figure}[h!]
\caption{Example of First Return - $M_{n}(x)$}
\end{figure}
\end{center}
\pagebreak
The path $M_{m,n}$, which denotes the Motzkin path with both $m$-colored down steps and $n$-colored level steps, will be represented by:

\begin{align}
    M_{m,n} = 1+nx\ M_{m,n}+mx^2\ M^2_{m,n}
\end{align}

\begin{center}
\begin{tikzpicture}[ultra thick]
\draw[->] [color=gray] (0,0)--(0,4);
\draw[->] [color=gray] (0,0)--(10,0);
\draw[->] (5,1)--(6,0);
\draw [thick, nearly transparent, dashed] (1,1) -- (5,1);
\draw[color=teal] (1,1) sin (1.5,2.1) cos (2,1.3) sin (2.5,2.5) cos (3,1.9) sin (3.5,3.5) cos (4,1.85) sin (4.5,2.7) cos (5,1);
\draw [->] (0,0)--(1,1) node[midway,above]{\large{$x$} \ \ };
\draw[->] (5,1)--(6,0);
\draw[->, red] (0,0.09)--(1,0.09);
\draw[->, green] (0,0.025)--(1,0.025);
\draw[->, yellow] (0,-0.05)--(1,-0.05);
\draw[->, blue] (0,-0.12)--(1,-0.12);
\draw[-] (0.5,-0.15)--(0.5,-0.15) node[midway, below]{\small{$nx$}};
\draw[->, red] (5.1,1)--(6.1,0);
\draw[->, green] (5,1)--(6,0);
\draw[->, yellow] (4.9,1)--(5.9,0);
\draw[->, blue] (4.8,1)--(5.8,0);
\draw[-] (5.65,0.55)--(5.65,0.55) node[midway, above]{\small{$mx$}};
\draw[color=teal] (6,0) sin (6.5,1.5) cos (7,0.25) sin (7.5,2.85) cos (8,0.25) sin (8.5,1.85) cos (9,0.4) sin (9.5,1.5) cos (10,0);
\draw[double, dotted, color=green] (1,0) sin (1.5,1) cos (2,0.2) sin (2.5,1.4) cos (3,0.5) sin (3.5,1) cos (4,0.1) sin (4.5,0.9) cos (5,0.3) sin (5.5,1.75) cos (6,1.15) sin (6.5,2.7) cos (7,.35) sin (7.5,1.35) cos (8,0.5) sin (8.5,0.1) cos (9,0.6) sin (9.5,2.5) cos (10,0);
\draw[-] (2.2,2.6)--(2.2,2.6) node[midway, above]{\large{$M_{m,n}(x)$}};
\draw[-] (7.4,3)--(7.4,3) node[midway, above]{\large{$M_{m,n}(x)$}};
\draw[-] (3.5,1)--(3.5,1) node[midway, above]{\large{$M_{m,n}(x)$}};

\end{tikzpicture}
\begin{figure}[h!]
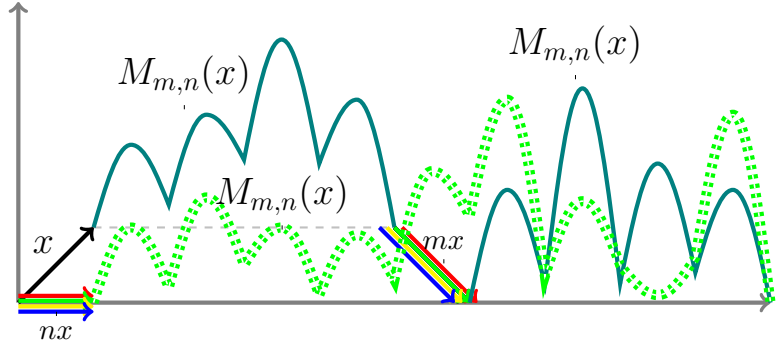

\caption{Example of First Return - $M_{m,n}(x)$}
\end{figure}
\end{center}

After using the quadratic formula, the respective generating functions yielded for $M_m$, $M_n$, and $M_{m,n}$ are:

\begin{align}
    M_m = \frac{1-x-\sqrt{1-2x+x^2(1-4m)}}{2mx^2}
\end{align}

\begin{align}
    M_n = \frac{1-nx-\sqrt{1-2nx+x^2(n^2-4)}}{2x^2}
\end{align}

\begin{align}
    M_{m,n} = \frac{1-nx-\sqrt{1-2nx+x^2(n^2-4m)}}{2mx^2}
\end{align}

The generating function for the Motzkin path is stated below for comparison purposes:
\begin{align}
    M(x)=\frac{1-x-\sqrt{1-2x-3x^2}}{2x^2} \tag{\ref{eq:4}}
\end{align}

\section{Conclusion}
Since we were unable to find the generating functions for the multicolored small Schr\"{o}der paths $s_n$ and $s_{m,n}$, one of us is investigating further into the problem.

\section{Acknowledgements}
This research is made possible by the generous support of the National Security Agency (NSA), Mathematical Association of America (MAA), National Science Foundation (NSF) grant DMS-1560332 administered though the American Statistical Association (ASA), Delta Kappa Gamma Educational Foundation, and Morgan State University.

\end{document}